\documentclass[12pt]{article}
\usepackage{amssymb}
\newtheorem{defn}{Definition}[section]
\newtheorem{lemma}[defn]{Lemma}
\newtheorem{prop}[defn]{Proposition}
\newtheorem{thm}[defn]{Theorem}

\newtheorem{rmk}[defn]{Remark}
\newtheorem{ex}[defn]{Example}

\def\C{{\mathbb C}}

\def\P{{\mathbb P}}

\def\II{{\cal I}}
\def\OO{{\cal O}}

\def\Bs{\mathop{\rm Bs}\nolimits}

\def\Hom{\mathop{\rm Hom}\nolimits}

\def\Pf{\mathop{\rm Pf}\nolimits}
\def\deg{\mathop{\rm deg}\nolimits}
\def\ker{\mathop{\rm ker}\nolimits}
\def\im{\mathop{\rm im}\nolimits}

\def\coker{\mathop{\rm coker\, }\nolimits}
\def\dim{\mathop{\rm dim\, }\nolimits}
\def\codim{\mathop{\rm codim\, }\nolimits}
\def\rk{\mathop{\rm rank\, }\nolimits}
\def\det{\mathop{\rm det}\nolimits}
\def\comp{{\scriptstyle \circ}}
\def\id{\mathop{\rm id}\nolimits}
\def\textwedge{{\textstyle\bigwedge}}

\def\ev{\mathop{\rm ev}\nolimits}

\def\mapright#1{\mathop{\vbox{\ialign{
                ##\crcr
    ${\scriptstyle\hfil\;\;#1\;\;\hfil}$\crcr
 \noalign{\kern-1pt\nointerlineskip}
    \rightarrowfill\crcr}}\;}}

\def\mapleft#1{\mathop{\vbox{\ialign{
                ##\crcr
    ${\scriptstyle\hfil\;\;#1\;\;\hfil}$\crcr
    \noalign{\kern-1pt\nointerlineskip}
    \leftarrowfill\crcr}}\;}}

\def\mapdown#1{\Big\downarrow
         \rlap{$\vcenter{\hbox{$\scriptstyle#1$}}$}}

\def\mapup#1{\Big\uparrow
         \rlap{$\vcenter{\hbox{$\scriptstyle#1$}}$}}

\def\proof{
  \noindent
  {\bf Proof:}
}

\def\endproof{
{\unskip\nobreak\hfil\penalty50\hskip2em\hbox{}\nobreak\hfill
          $\square$\bigbreak}
}

\begin{document}

\title{Non--vanishing for Koszul cohomology of curves}
\author{M. Aprodu, J. Nagel}
\date{}
\maketitle
\abstract{
We study the relationship between rank $p+2$ Koszul classes and rank 2 vector bundles on a smooth curve. We show that every rank $p+2$ Koszul class is obtained from a rank 2 vector bundle and give
an explicit nonvanishing theorem for Koszul classes arising in this way.
}

\section{Introduction}
\label{section: introduction}

Let $X$ be a smooth complex projective variety. The geometry of $X$ is
reflected in the behaviour of the Koszul cohomology groups $K_{p,q}(X,L)$
introduced by Green \cite{GreenJDG}, more specifically the
vanishing/nonvanishing of certain Koszul cohomology groups. The
fundamental result in this direction is the nonvanishing theorem of
Green--Lazarsfeld \cite{GL1}. This theorem states that if a line bundle
$L$ admits a decomposition $L=L_1\otimes L_2$ with $r_i = h^0(X,L_i)-1\ge 1$
($i=1,2$) then $K_{r_1+r_2-1,1}(X,L)\ne 0$. Voisin \cite[(1.1)]{VoisinLMS}
has given a different proof of this result under the hypothesis that $L_1$
and $L_2$ are globally generated.

\medskip
The aim of this note is to give a more geometric approach to this type of problems. The starting point is the following
construction due to Voisin. Given a rank two vector bundle $E$ on $X$ with determinant $L$,
Voisin \cite[(2.22)]{Voisin2} defined a homomorphism
$$
\varphi: S^p H^0(X,E)\otimes \textwedge^{p+2}H^0(X,E)\rightarrow \textwedge^p H^0(X,L)\otimes H^0(X,L).
$$
By \cite[Lemma 5]{Voisin2}, this homomorphism produces elements of $K_{p,1}(X,L)$. If we take
$E =L_1\oplus L_2$, we get back the classes constructed by Green and Lazarsfeld.
As one of the referees pointed out to us, Koh and Stillman \cite{KS} had generalised the Green--Lazarsfeld construction before from a different point of view.

\medskip
Recall that the {\em rank} of a Koszul class $\gamma\in K_{p,1}(X,L)$ is
the minimal dimension of a linear subspace $W\subset H^0(X,L)$ such that $\gamma$ is represented by an element in $\textwedge^p W\otimes H^0(X,L)$; cf. \cite[Definition 2.2]{Hans-Chris}. (Note that the subspace $W$ is uniquely determined if $p\ge 2$.) By definition, the Koszul classes constructed in this paper are of rank $p+2$ if the vector bundle $E$ is indecomposable.

\medskip
Section \ref{section: main} contains the main results of this paper.
We first give a necessary and sufficient condition for nonvanishing of Koszul classes on smooth curves obtained from rank 2 vector bundles (Theorem \ref{thm: main}). This result generalises the nonvanishing theorem of Green--Lazarsfeld in the case of curves. Our second main result, Theorem \ref{thm: Grassmann}, states that every rank $p+2$ Koszul class on a smooth curve comes from a rank two vector bundle. This theorem is a generalisation of \cite[Theorem 6.7]{Hans-Chris}.

\section{Preliminaries}
\label{section: prelim}

\subsection{The method of Voisin}
\label{subsec: Voisin}

Let $E$ be a rank two vector bundle on a smooth projective variety $X$ defined over an algebraically closed field $k$ of characteristic zero.
Write $L=\det E$ and $V = H^0(X,L)$, and let
$$
d:\textwedge^2 H^0(X,E)\rightarrow V
$$
be the determinant map.
Given $t\in H^0(X,E)$, define a linear map
$$
d_t:H^0(X,E)\rightarrow V
$$
by $d_t(u) = d(t\wedge u)$, and choose a subspace $U\subset H^0(X,E)$ with $U\cap\ker(d_t)=0$. Suppose that $\dim(U)=p+2$ with $p\ge 1$, and put $W = d_t(U)\cong U$. The restriction of $d$ to $\textwedge^2 U$ defines a map $\textwedge^2 U\rightarrow V$, which we can view as an element of
$$
\textwedge^2 U^{\vee}\otimes V\cong\textwedge^p U\otimes V.
$$
Let
$$
\gamma\in\textwedge^p W\otimes V\subset\textwedge^p V\otimes V
$$
be the image of this element under the map $d_t$.

\medskip

Following Voisin \cite[(2.22)]{Voisin2},
we prove that $\gamma$ defines a Koszul class in $K_{p,1}(X,L)$.
To this end, we make the previous construction explicit using coordinates.
If we choose a basis $\{e_1,\ldots,e_{p+3}\}$ of $\langle t\rangle \oplus U\subset H^0(X,E)$ such that $e_1 = t$, we have
\begin{eqnarray}
\label{eqn: gamma}
\gamma & = & \sum_{i<j}(-1)^{i+j}d(t\wedge e_2)\wedge\ldots\wedge\widehat{d(t\wedge e_i)}\wedge\ldots \\
& & \nonumber \ldots\wedge\widehat{d(t\wedge e_j)}\wedge\ldots\wedge d(t\wedge e_{p+3})\otimes d(e_i\wedge e_j).
\end{eqnarray}
As in \cite{Voisin2} one shows that the image of the $\gamma$ by the Koszul differential
$$
\delta:\textwedge^p V\otimes H^0(X,L)\rightarrow\textwedge^{p-1}V\otimes S^2H^0(X,L)
$$
equals
\begin{eqnarray}
\label{eqn:delta(gamma)}
\sum_{i<j<k} (-1)^{i+j+k} d(t\wedge e_2)\wedge\ldots \widehat{d(t\wedge e_i)}\ldots\widehat{d(t\wedge e_j)}\ldots\widehat{d(t\wedge e_k)}\ldots\wedge{d(t\wedge e_{p+3})} \\
\nonumber
\otimes
\{d(t\wedge e_i)d(e_j\wedge e_k)-d(t\wedge e_j)d(e_i\wedge e_k)+d(t\wedge e_k)d(e_i\wedge e_j)\}.
\end{eqnarray}

\begin{lemma}[Voisin]
\label{lemma: Voisin}
Given four elements $w_1$, $w_2$, $w_3$, $w\in H^0(X,E)$ we have the relation
$$
d(w\wedge w_1)d(w_2\wedge w_3) - d(w\wedge w_2)d(w_1\wedge w_3)+
d(w\wedge w_3)d(w_1\wedge w_2) = 0
$$
in $H^0(X,L^2)$.
\end{lemma}

\proof
See \cite[Lemma 5]{Voisin2}.
\endproof

\medskip
The previous lemma shows that $\gamma$ belongs to the kernel of the Koszul differential
$$
\delta_X:\textwedge^{p}V\otimes H^0(X,L)\rightarrow\textwedge^{p-1}V\otimes H^0(X,L^2).
$$
Hence $\gamma$ defines a Koszul class $[\gamma]\in K_{p,1}(X,L,W)\subseteq K_{p,1}(X,L)$.
Clearly the given class only depends on $t$ and $W$; we write $[\gamma]=\gamma(W,t)$.

\subsection{The method of Green--Lazarsfeld}
\label{subsec: GL}

\medskip
Let $L_1$, $L_2$ be two line bundles on a smooth projective variety $X$ such that $r_i = h^0(X,L_i)-1\ge 1$ ($i=1$, 2). Write $L_i = M_i+F_i$ with $M_i$ the mobile part and $F_i$ the fixed part. Let $B$ be the divisorial part of $F_1\cap F_2$. It is possible to choose $s_i\in H^0(X,L_i)$ such that $V(s_1,s_2) = B\cup Z$ with $\codim(Z)\ge 2$. Set $L = L_1\otimes L_2$, and put $t=(s_1,s_2)\in H^0(X,L_1\oplus L_2)$, $W = \im(d_t)\subset H^0(X,L(-B))$. By construction
$h^0(X,\OO_X(B))=1$, hence $\dim W = r_1+r_2+1$. By the previous discussion, we obtain a Koszul class $\gamma(W,t)\in K_{r_1+r_2-1,1}(X,L)$. We call such classes {\em Green--Lazarsfeld classes}.

\medskip
Note that the rank of a Green--Lazarsfeld class is either $p+1$ or $p+2$. Classes of rank $p+1$ are of scrollar type; see e.g. \cite{Schreyer} or \cite[Corollary 5.2]{Hans-Chris}.

\begin{defn}
\label{def: GL}
Given a nonnegative integer $k\ge 0$, let $K_{k,1}(X,L)_{\rm GL}\subseteq K_{k,1}(X,L)$ be the subspace generated by Green--Lazarsfeld classes for all decompositions $L = L_1\otimes L_2$
with $k= r_1+r_2-1$, ($r_1\ge 1$, $r_2\ge 1$).
\end{defn}

\subsection{The method of Koh--Stillman}
\label{subsec: KS}

\medskip
Voisin's method produces syzygies of rank $\le p+2$. As we have seen in the previous subsection, rank $p+1$ syzygies are Green--Lazarsfeld syzygies of scrollar type. Rank $p+2$ syzygies can be obtained in the following way. Suppose that $L$ is a globally generated line bundle on a projective variety $X$, and let  $[\gamma]\in K_{p,1}(X,L)$ be a nonzero class represented by an element $\gamma\in\textwedge^p W\otimes V$ with $\dim W = p+2$. We view $\gamma$ as an element in $\textwedge^2 W^{\vee}\otimes V\cong\Hom(\textwedge^2 W,V)$. Following \cite[Proof of Theorem 6.1]{Hans-Chris} we consider the map
$$
\gamma^{\prime}: \textwedge^2(\C\oplus W)=W\oplus\textwedge^2 W\rightarrow V
$$
defined by taking the direct sum of $\gamma$ and the inclusion $W\hookrightarrow V$. If we choose a generator $e_1$ for the first summand and a basis $\{e_2,\ldots,e_{p+3}\}$ for $W$, we obtain a skew--symmetric $(p+3)\times(p+3)$ matrix $A$ by setting
$$
a_{ij} = \gamma^{\prime}(e_i\wedge e_j).
$$
By construction, the inclusion $W\rightarrow V$ corresponds to the map $\gamma^{\prime}(e_1\wedge -)$. This allows us to identify $a_{1j}$ and $e_j$, $2\le j\le p+3$.
Let $\alpha$ be the image of $\gamma$ under the Koszul differential
$$
\delta:\textwedge^p V\otimes V\rightarrow\textwedge^{p-1}V\otimes S^2 V.
$$
Writing this out, we obtain
\begin{equation}
\label{eq: KS} \alpha = \sum_{i<j<k}
(-1)^{i+j+k}a_{12}\wedge\ldots\widehat{a_{1,i}}\ldots\widehat{a_{1,j}}
\ldots\widehat{a_{1,k}}\ldots \wedge a_{1,p+3}\otimes\Pf_{1ijk}(A).
\end{equation}
As the elements $\{a_{12},\ldots a_{1,p+3}\} = \{e_2,\ldots,e_{p+3}\}$ are linearly independent, this expression is nonzero if and only if at least one of the Pfaffians $\Pf_{1ijk}(A)$ is nonzero. Furthermore, since $\alpha$ maps to zero in $\textwedge^{p-1}V\otimes H^0(X,L^2)$ the Pfaffians
$\Pf_{1ijk}(A)$ have to vanish on the image of $X$.

\medskip
The preceding discussion shows that every rank $p+2$ syzygy arises from a skew--symmetric $(p+3)\times (p+3)$ matrix $A$ such that
\begin{enumerate}
\item[(i)] the elements $\{a_{12},\ldots a_{1,p+3}\}$ are linearly independent;
\item[(ii)] there exists a nonzero Pfaffian $\Pf_{1ijk}(A)$;
\item[(iii)] the Pfaffians $\Pf_{1ijk}(A)$ vanish on the image of $X$ in $\mathbb{P}(V^\vee)$.
\end{enumerate}
This is exactly the method used by Koh and Stillman to produce syzygies; see \cite[Lemma 1.3]{KS}.

\begin{rmk}
 \label{rmk: Eisenbud}
 {\rm In the geometric setting of subsection \ref{subsec: Voisin},
 let $Y$ be the image of $X$ in $\mathbb{P}(V^\vee)$. The expression (\ref{eqn:delta(gamma)}) shows that the canonical isomorphism
 $$
 K_{p,1}(X,L)\cong K_{p-1,2}(\P^r,\II_Y,\OO_{\P}(1))
 $$
 maps the class $\gamma(W,t)$ to the element $\alpha$ defined in (\ref{eq: KS}). Moreover, if $d$ does not vanish on decomposable elements then $\gamma(W,t)\ne 0$.
 Indeed, this condition is satisfied if and only if the matrix $A$ has no generalised zero;
 cf. \cite[Definition (1.1)]{KS}. One then applies [loc. cit., Remark p. 122].
 }
\end{rmk}

\medskip

\section{Main results}
\label{section: main}

\medskip
\begin{thm}
\label{thm: main}
Let $X$ be a smooth curve, let $L$ be a base--point free line bundle on $X$ and let $W\subset H^0(X,L)$ be a linear subspace. Put $B = \Bs(W)$, and let $t$ be a section of $H^0(X,\OO_X(B))$ vanishing on $B$. Consider an extension
\begin{equation}
\label{eqn:extension}
0\rightarrow \OO_X(B)\rightarrow E\rightarrow L(-B)\rightarrow 0
\end{equation}
such that
$$
W\subset(\ker H^0(X,L(-B))\mapright{\delta} H^1(X,\OO_X(B))).
$$
Then the Koszul class $\gamma(W,t)$ defined in section \ref{subsec: Voisin} is nonzero
is and only if the extension (\ref{eqn:extension}) is non-split.
\end{thm}

\proof
The proof proceeds in several steps. We use the notation of section \ref{subsec: Voisin}.

\medskip\noindent{\bf Step 1.}
Suppose that the extension (\ref{eqn:extension}) splits. In this case, one readily verifies that $d$ vanishes identically on $\textwedge^2 U$. The formula (\ref{eqn: gamma}) then shows that $\gamma(W,t) = 0$.

\medskip\noindent{\bf Step 2.}
If $\gamma(W,t) = 0$ there exists a linear map $h:U\rightarrow \C$ such that
\begin{equation}
\label{eq: linear map}
d(u_1\wedge u_2) = h(u_2)d_t(u_1)-h(u_1)d_t(u_2)
\end{equation}
for all $u_1,u_2\in U$.

\medskip
Indeed, suppose that there exists a nonzero element $\tilde\gamma\in\textwedge^{p+1}W\cong W^{\vee}$ such that $\gamma$ is the image of $\tilde\gamma$ under the Koszul differential. Then $\gamma$ coincides with the composition of maps
$$
\textwedge^2 W\mapright{\delta}W\otimes W\mapright{\tilde\gamma\otimes\id}
W\hookrightarrow V.
$$
Since
\begin{eqnarray*}
d(u_1\wedge u_2) & = & \gamma(d_t(u_1)\wedge d_t(u_2)) \\
& = & \tilde\gamma(d_t(u_2))d_t(u_1)-\tilde\gamma(d_t(u_1))d_t(u_2)),
\end{eqnarray*}
condition (\ref{eq: linear map}) is satisfied with
$h=\tilde\gamma\comp d_t:U\rightarrow\C$.

\medskip\noindent{\bf Step 3.}
Let $u_1$, $u_2\in U$ be two sections such that $d_t(u_1)$ and $d_t(u_2)$ generate $L(-B)$. If $d(u_1\wedge u_2) = 0$, the extension (\ref{eqn:extension}) splits.

\medskip
To prove this assertion, put $s_i = d_t(u_i)$ ($i=1,2$) and consider the commutative diagram
$$
\begin{array}{ccccccccc}
0  & \rightarrow & \OO_X(B) & \mapright{}  & E & \mapright{} & L(-B) & \rightarrow & 0\\
& & & & \mapup{\ev_1} & & \mapup{\ev_2} & & \\
&& 0 &  \rightarrow & \langle u_1\, , u_2 \rangle\otimes \OO_X & \stackrel{\sim}{\rightarrow}
& \langle s_1\, ,s_2\rangle \otimes \OO_X & \rightarrow & 0.
\end{array}
$$
Put $M=\ker(\ev_1)$, and note that $\ker(\ev_2)\cong L^{-1}(B)$ since $\ev_2$
is surjective. By the Snake Lemma we obtain an exact sequence
$$
0\rightarrow M \rightarrow L^{-1}(B) \rightarrow\OO_X(B)
\rightarrow \coker(\ev_1) \rightarrow 0.
$$
Note that
$$
d(u_1\wedge u_2)=0 \Longleftrightarrow \rk
\im (\langle u_1\, ,u_2\rangle \otimes \OO_X
\rightarrow E)=1\Longleftrightarrow \rk M = 1
$$
where the first equivalence follows from \cite[p. 380]{Voisin1}.
If $d(u_1\wedge u_2) = 0$ the above exact sequence shows that $M\cong L^{-1}(B)$, hence the isomorphism $\langle u_1\, ,u_2 \rangle\otimes \OO_X\stackrel{\sim}{\rightarrow}
\langle s_1\, ,s_2\rangle \otimes \OO_X$ induces an isomorphism
$\im(\ev_1)\cong L(-B)$. The inverse
of this isomorphism provides a splitting of the extension (\ref{eqn:extension}).

\medskip\noindent{\bf Step 4.}
Suppose that $\gamma(W,t) = 0$. Then there exists a linear map $h:U\rightarrow\C$ as in Step 2. Consider the morphism
$$
\pi:X\rightarrow\P(W^{\vee})
$$
defined by the base--point free linear system $W\subset H^0(X,L(-B))$, and choose a linear subspace $\Lambda\subset\P(W^{\vee})$ of codimension two such that $\Lambda\cap\pi(X) = \emptyset$. The hyperplane $\ker(h)\subset W$ corresponds to a point $p\in\P(W^{\vee})$. Put $H_1 = \langle\Lambda, p\rangle$ and choose a hyperplane $H_2\subset\P(W^{\vee})$ containing $\Lambda$ such that $p\notin H_2$. Let $u_1$, $u_2$ be the sections corresponding to $H_1$, $H_2$. Then $d_t(u_1)$ and $d_t(u_2)$ generate $L(-B)$ and $u_1\in\ker(h)$, $u_2\notin\ker(h)$.
Equation (\ref{eq: linear map}) yields the identity
$$
d(u_1\wedge u_2) = h(u_2)d_t(u_1).
$$
Rewriting this identity, we obtain $d(u_1\wedge(u_2 + h(u_2)t)) = 0$. Since the pair $\{d_t(u_1),d_t(u_2+h(u_2)t)\} = \{d_t(u_1),d_t(u_2)\}$ generates $L(-B)$, Step 3 implies that the extension (\ref{eqn:extension}) splits.
\endproof

\medskip
\begin{rmk}
\label{rmk: refinement}
{\rm
In the statement of Theorem \ref{thm: main} it is not necessary to suppose that $L$ is globally generated, since $K_{p,1}(X,L(-\Bs(L)))\cong K_{p,1}(X,L)$.
}
\end{rmk}

\medskip
Theorem \ref{thm: main} yields a short, geometric proof of
the Green--Lazarsfeld nonvanishing theorem for curves.
\begin{thm} {\bf (Green--Lazarsfeld)}
\label{cor: GL}
Let $X$ be a smooth curve, and let $L$ be a line bundle on $X$ that admits
a decomposition $L=L_1\otimes L_2$ with $r_i=\dim|L_i|\ge 1$ for $i=1,2$.
Then $K_{r_1+r_2-1,1}(X,L)\ne 0$.
\end{thm}

\proof
We define $s_1$, $s_2$, $t$, $W$, $B$ and $\gamma(W,t)$ as in section \ref{subsec: GL}. Let $C$ be the base locus of $W$, seen as a subspace of $H^0(X,L(-B))$. We prove that $\gamma(W,t)\ne 0$.
Suppose that $\gamma(W,t)=0$. Consider the extension
$$
0\rightarrow\OO_X(B)\rightarrow L_1\oplus L_2\rightarrow L(-B)\rightarrow 0.
$$
Pulling back this extension along the injective homomorphism $L(-B-C)\rightarrow L(-B)$, we obtain an induced extension
$$
0\rightarrow \OO_X(B)\rightarrow E\rightarrow L(-B-C)\rightarrow 0.
$$
Applying Theorem \ref{thm: main} to the line bundle $L(-C)$, we find that this extension splits. Hence there exists an injective homomorphism
$$
\OO_X(B)\oplus L(-B-C)\rightarrow L_1\oplus L_2.
$$
In particular there exists $i\in\{1,2\}$ such that $\Hom(L(-B-C),L_i)\ne 0$.
This implies that
$$
r_i+1 = h^0(X,L_i)\ge h^0(X,L(-B-C))\ge\dim W = r_1+r_2+1,
$$
and this is impossible since $r_1\ge 1$ and $r_2\ge 1$.
\endproof

\begin{thm}
\label{thm: Grassmann}
Let $X$ be a smooth curve, and let $\alpha\ne 0\in K_{p,1}(X,L)$ be a Koszul class of rank $p+2$ represented by an element of $\textwedge^p W\otimes H^0(X,L)$ with $\dim W = p+2$. There exist a rank 2 vector bundle $E$ on $X$ and a section $t\in H^0(X,E)$ such that $\alpha = \gamma(W,t)$.
\end{thm}

\proof
Put $T = \C\oplus W$, and choose a basis $\{e_1,\ldots,e_{p+3}\}$ of $T$ such that $t=e_1$ is the generator of the first summand. Writing $z_{ij} = e_i\wedge e_j$, we obtain a skew--symmetric matrix $Z = (z_{ij})$ and
coordinates $(z_{ij})_{1\le i<j\le p+3}$ on $\P(\textwedge^2 T^{\vee})$.
Consider the Grassmannian $G = G(2,T)$ of 2--dimensional quotients of $T$. The ideal of $G$ under the Pl\"ucker embedding $G\subset\P(\textwedge^2 T^{\vee})$ is generated by the $4\times 4$ Pfaffians $\Pf_{ijkl}(Z)$ of the matrix $Z$. Taking exterior powers in the exact sequence
$$
0\rightarrow \langle t\rangle\rightarrow T\rightarrow W\rightarrow 0
$$
we obtain an exact sequence
$$
0\rightarrow\langle t\rangle\otimes W\rightarrow\textwedge^2 T\rightarrow\textwedge^2 W\rightarrow 0.
$$
The linear subspace $\P(\textwedge^2 W^{\vee})\subset\P\left(\textwedge^2 T^{\vee}\right)$ is defined by the vanishing of the linear forms $z_{1j}$, $j=2,\ldots,p+3$. A straightforward computation then shows that the ideal of the union
$$
G(2,T)\cup\P(\textwedge^2 W^{\vee})\subset\P(\textwedge^2 T^{\vee})
$$
is generated by the Pfaffians $\Pf_{1ijk}(Z)$.
The tautological exact sequence
$$
0\rightarrow S\rightarrow T\otimes\OO_G\rightarrow Q\rightarrow 0
$$
induces an isomorphism $T\cong H^0(G,Q)$. Under this isomorphism, we have $G(2,W) = V(t)$.
\medskip
As in section \ref{subsec: KS} we associate to the Koszul class $\alpha$ a matrix $A = (a_{ij})$ of linear forms $A = (a_{ij})$ such that
\begin{enumerate}
\item[(a)] The linear forms in the first row of $A$ span $W$;
\item[(b)] There exists a nonzero $4\times 4$ Pfaffian of $A$
involving the first row and column;
\item[(c)] The $4\times 4$ Pfaffians involving the first row and column of $A$ vanish on the image of $X$ in $\mathbb{P} H^0(X,L)^\vee$.
\end{enumerate}
Let $C$ be the base locus of the image of $A$. Replacing $L$ by $L(-C)$ if necessary ($W$ is obviously contained in the image of $A$) we can suppose that $C$ is empty, hence the matrix $A$ defines a morphism
$$
\psi:X\rightarrow\P(\textwedge^2 T^{\vee}).
$$
Condition (c) implies that the image $Y = \psi(X)$ is contained in the union $G(2,T)\cup\P(\textwedge^2 W^{\vee})$, and condition (a) shows that $Y$ is not contained in $\P(\textwedge^2 W^{\vee})$. As $Y$ is irreducible, this implies that $Y$ is contained in $G(2,T)$.

\medskip
Put $E = \psi^*Q$. Twisting the exact sequence
$$
0\rightarrow\II_Y\rightarrow\OO_G\rightarrow\psi_*\OO_X\rightarrow 0
$$
by the universal quotient bundle $Q$ and taking global sections, we obtain an exact sequence
$$
0\rightarrow H^0(G,Q\otimes\II_Y)\rightarrow H^0(G,Q)\mapright{\psi^*} H^0(G,\psi_*\OO_X\otimes Q)\cong H^0(X,E).
$$
Condition (a) implies that $Y$ is not contained in
$G(2,W) = G(2,T)\cap\P(\textwedge^2 W^{\vee})$, hence $t$ does not vanish identically on $X$ and defines a global section of $E$. The zero locus of this section is given by the equations $a_{12}=\dots=a_{1,p+3}=0$, hence it coincides with $B$. Consequently the line bundle $E$
is given by an extension
\begin{equation}
\label{eqn: Grass ext}
0\to \OO_X(B)\to E\to L(-B)\to 0.
\end{equation}

\medskip
Consider the commutative diagram
$$
\begin{array}{ccc}
0 & & 0 \\
\mapdown{} & & \mapdown{} \\
H^0(G,\OO_G) & \mapright{.t} & H^0(X,\OO_X(B)) \\
\mapdown{\wedge t} & & \mapdown{\wedge t} \\
H^0(G,Q) & \mapright{\psi^*} & H^0(X,E) \\
\mapdown{} & & \mapdown{d_t} \\
W & \mapright{i} & H^0(X,L(-B)).
\end{array}
$$
Note that $\ker i = W\cap H^0(G,\OO_G(1)\otimes\II_Y) = 0$ by condition (a). As the map $H^0(G,Q)\rightarrow W$ is surjective, we find that $W$ is contained in the image of the map $d_t:H^0(X,E)\to H^0(X,L(-B))$. Hence the condition of Theorem \ref{thm: main} is satisfied.
By condition (b) we have $\gamma(W,t)\ne 0$. Hence the extension (\ref{eqn: Grass ext}) does not split by Theorem \ref{thm: main}.
\endproof

\begin{rmk}
{\rm
The union $G(2,T)\cup\P(\textwedge^2 W^{\vee})$ is a generic syzygy scheme; see \cite[Theorem 6.1]{Hans-Chris}. In [loc. cit., Theorem 6.7] it was shown that a rank $p+2$ syzygy gives rise to a rank 2 vector bundle if $L$ is very ample and the ideal of $X$ is generated by quadrics.
}
\end{rmk}

\medskip
The condition of Theorem \ref{thm: main} can be reinterpreted in terms of surjectivity of a natural multiplication map.
\begin{prop}
\label{prop: extension}
Let $X$ be a smooth curve, and let $W\subset H^0(X,L)$ be a linear subspace. We put $B = \Bs(W)$ and view $W$ as a base--point free linear subspace of $H^0(X,L(-B))$. Let
$$
\mu:W\otimes H^0(X,K_X(-B))\rightarrow H^0(K_X\otimes L(-2B))
$$
be the multiplication map. The following conditions are equivalent.
\begin{enumerate}
\item[(i)] The map $\mu$ is not surjective;
\item[(ii)] There exists a non-split
extension
$$
0\rightarrow \OO_X(B)\rightarrow E\rightarrow L(-B)\rightarrow 0
$$
such that $W$ is contained in the kernel of the map
$\delta:H^0(X,L(-B))\rightarrow H^1(X,\OO_X(B))$.
\end{enumerate}
\end{prop}

\proof
We first show that (i) implies (ii). Since $\mu$ is not surjective, there exists a
hyperplane $H\subset H^0(X,K_X\otimes L(-B))$ that contains $\im(\mu)$. Let $\eta$ be a linear functional defining $H$. Put $0\ne\xi=\eta^\vee\in H^1(X,L^{-1}(B))$, and let
$$
0\rightarrow \OO_X(B)\rightarrow E\rightarrow L(-B)\rightarrow 0
$$
be the corresponding non-split extension.
Given $w\in W$ and $v\in H^0(X,K_X(-B))$, the formula
\begin{equation}
\label{eq: delta}
\delta(w)(v)=(\eta\circ\mu)(w\otimes v)
\end{equation}
shows that $W$ is contained in the kernel of $\delta$.

\medskip
For the converse, note that formula (\ref{eq: delta}) implies that
$\eta|_{\im\mu}\equiv 0$.
\endproof

\begin{rmk}
{\rm
If $B$ is a fixed divisor, the result of the previous Proposition follows from Green's duality theorem \cite[Corollary (2.c.10)]{GreenJDG}. Indeed,
\begin{equation}
\label{eqn: duality}
\coker \mu \cong K_{0,1}(X,K_X(-B),L(-B),W)\cong
K_{p,1}(X,B,L(-B),W)^\vee
\end{equation}
and since $h^0(X,\OO_X(B))=1$ we have an injection
$$
K_{p,1}(X,B,L(-B),W)\hookrightarrow K_{p,1}(X,L).
$$
}
\end{rmk}

\medskip
Theorem \ref{thm: Grassmann} shows that Voisin's method may produce nontrivial Koszul classes that are not contained in the space $K_{p,1}(X,L)_{\rm GL}$ spanned by Green--Lazarsfeld classes.
\begin{ex}
{\rm
By \cite[Theorem 3.6 and Theorem 4.3]{ELMS} there exists a smooth curve of genus 14 and Clifford index 5 whose Clifford index is computed by a unique line bundle $L$ such that $L^2 = K_X$. The line bundle $L$ embeds $X$ in $\P^4$ as a projectively normal curve of degree 13 which
is not contained in any quadric of rank $\le 4$, and the ideal of $X$ is generated by the $4\times 4$ Pfaffians of a skew--symmetric matrix $(a_{ij})_{1\le i,j\le 5}$ with
$$
\deg(a_{ij})=\left\{
\begin{array}{l}
 2\mbox{ if } i=1    \mbox{ or }  j=1\\
 1\mbox{ if } i\ge 2 \mbox{ and } j\ge 2
\end{array}
\right.
$$
such that the quadric $Q=a_{23}a_{45}-a_{24}a_{35}+a_{25}a_{34}$
has rank 5.

\medskip
By [loc.cit.] the group $K_{1,1}(X,L)$ is generated by $[Q]$, hence $I_X$ contains no quadrics of rank $\le 4$. If $K_{1,1}(X,L)$ contains a Green--Lazarsfeld class this class would be of scrollar type, since it necessarily comes from two pencils $|L_1|$, $|L_2|$. This is impossible, since classes of scrollar type give rise to quadrics of rank $\le 4$.

\medskip
The Koszul class $[Q]\in K_{1,1}(X,L)$ has rank 3, since it is represented by the linear subspace $W = \langle a_{23},a_{24},a_{25}\rangle$. Hence $[Q]$ comes from Voisin's method by Theorem \ref{thm: Grassmann}.
}
\end{ex}

\begin{rmk}
{\rm
A more geometric description of a subspace $W$ representing $[Q]$ is the following.
A smooth intersection of the quadric
$V(Q)\subset\mathbb{P}H^0(X,L)^{\vee}$ with one of the cubic Pfaffians
is a $K3$ surface in $\mathbb{P}H^0(X,L)^{\vee}$ containing a line $\ell$ which is disjoint from $X$ by \cite[Prop. 4.1]{ELMS}. The line $\ell$ corresponds to a $3$-dimensional linear subspace $W\subset H^0(X,L)$, which
is base-point-free since $\ell$ does not meet $X$.
}
\end{rmk}

\medskip
One could ask whether the syzygies constructed in section \ref{subsec: Voisin} span $K_{p,1}(X,L)$. In principle it may be possible to obtain higher rank syzygies as linear combinations of rank $p+2$ syzygies. However, if $K_{p,1}(X,L)$ is spanned by a single syzygy of rank $\ge p+3$ this is not possible.
\begin{ex}[Eusen--Schreyer]
\label{ex: ES}
{\rm
Eusen and Schreyer  \cite[Theorem 1.7 (a)]{ES} have constructed a smooth curve $X\subset\P^5$ of genus 7 and Clifford index 3 embedded by the linear system $|K_X(-x)|$ such that $K_{2,1}(X,K_X(-x))\cong\C$ is spanned by a syzygy $s_0$. The explicit expression for $s_0$ given on p.8 of [loc. cit.]
shows that $s_0$ is a rank 5 syzygy. Hence $s_0$ cannot be obtained by the Green--Lazarsfeld construction or the method of section \ref{subsec: Voisin}.
}
\end{ex}

\noindent {\bf Acknowledgements.} The first named author was
partially supported by a Humboldt Research Fellowship, and by the
ANCS contract 2-CEx 06-11-20/2006. We would like to thank Universit\'e
Grenoble 1, I.H.E.S., Universit\"at Bayreuth and Universit\'e Lille
1 for hospitality during the first stage of this work.

\medskip
We thank the referees for several comments that helped to improve the presentation of the paper, and for pointing out an error in the previous version of this paper.

\noindent
\\\\
M. Aprodu, Romanian Academy, Institute of Mathematics "Simion
Stoilow", P.O.Box 1-764, RO-014700, Bucharest, ROMANIA, email:
Marian.Aprodu\char64 imar.ro
\\\\
\c Scoala Normal\u a Superioar\u a--Bucure\c sti
21, Calea Grivitei Str.
010702-Bucharest, Sector 1
ROMANIA
\\\\
J. Nagel, Universit\'e Lille 1,
Math\'ematiques -- B\^at. M2,
F-59655 Villeneuve d'Ascq Cedex,
FRANCE,
email: nagel\char64 math.univ-lille1.fr

\end{document}